\input amstex
\documentstyle{amsppt}
\mag=1200

\widestnumber\key{FrolTa}

\define\al{\alpha}
\define\be{\beta}

\define\e{\varepsilon}

\define\ov{\overline}

\define \Tr{\operatorname{Tr}}

\define\0{\bold0}

\define\BZ{\Bbb Z}
\define\BQ{\Bbb Q}
\define\BR{\Bbb R}
\define\BT{\Bbb T^m}
\define\BN{\Bbb N}

\define\A{\Cal A}

\define \({\left(}
\define \){\right)}

\redefine\t{\bold t}

\redefine\P{\Cal P}
\redefine\H{\Cal H}

\redefine\phi{\varphi}

\topmatter

\author Nikita Sidorov \endauthor

\title An arithmetic group associated with a Pisot unit and its
symbolic-dynamical representation
\endtitle

\rightheadtext{An arithmetic group associated with a Pisot unit}

\abstract To a given Pisot unit $\beta$ we associate a finite abelian group
whose size appears to be equal to the discriminant of $\beta$. We call
it the Pisot group and find its representation in the two-sided
$\beta$-compactum in the case of $\beta$ satisfying the relation
$Fin(\beta)=\Bbb Z[\beta]\cap[0,1)$. As a motivation for the definition,
we show that the Pisot group is the kernel of some important arithmetic
coding of the toral automorphism given by the companion matrix naturally
associated with $\beta$.
\endabstract


\address Department of Mathematics, UMIST, P.~O.~Box 88, Manchester M60 1QD,
United Kingdom
\endaddress

\email Nikita.A.Sidorov$\@$umist.ac.uk
\endemail

\date December 23, 1999 \enddate

\endtopmatter

\document

\head 1. The definition of the Pisot group and its basic properties \endhead

Let\footnote""{Supported by the EPSRC grant ref.~no.~GRL98923.}
$\be>1$ be a {\it Pisot number}, i.e. an algebraic integer whose
conjugates have the moduli strictly less than 1. Let the characteristic
polynomial of $\be$ be $g(x)=x^m-k_1x^{m-1}-k_2x^{m-2}-\dots-k_m$.
We assume $\be$ to be a {\it unit}, i.e. $k_m=\pm1$.

We recall that since $\be$ is a Pisot number, $\|\be^n\|\to0$ as
$n\to+\infty$, where $\|y\|=\min\,\{|y-k|,\ k\in\BZ\}$. Let
$$
\P_\be=\{\xi: \|\xi\be^n\|\to0,\,\, n\to+\infty\}.
$$

Let us first establish some auxiliary facts. It is well-known that
$\xi\in\BQ(\be)$ (see, e.g., \cite{Cas}).
Let $\Tr(\xi)$ denotes the trace of $\xi$, i.e. the sum of $\xi$
and all its conjugates.
\proclaim{Lemma 1.1} The set $\P_\be$ is a commutative group under addition containing
$\BZ[\be]$. It can be characterized as follows:
$$
\P_\be=\{\xi\in\BQ(\be):\Tr(a\xi)\in\BZ\ \text{\rm{for any}}\ a\in\BZ[\be]\}.
\tag1.1
$$
\endproclaim
\demo{Proof} The first claim is a consequence of the inequalities
$\|\xi_1\pm\xi_2\|\le\|\xi_1\|+\|\xi_2\|$ and of the fact that
$\|\be^n\|\to0$. To prove (1.1), we observe that
the classical theorem due to Pisot and Vijayaraghavan says that
for any Pisot $\be$, there exists $k_0\in\BN$ such that
$\P_\be=\{\xi\in\BQ(\be):\Tr(\be^k\xi)\in\BZ,\ k\ge k_0\}$ (see, e.g.,
\cite{Cas}), whence
follows (1.1), as $\be$ is a unit, and $\Tr(\zeta)\in\BZ$ implies
$\Tr(\be^{-1}\zeta)\in\BZ$.\qed
\enddemo

Thus, if we regard $\BZ[\be]$ as a lattice over $\BZ$, then by (1.1) and
the definition, $\P_\be$ is the dual lattice for $\BZ[\be]$ (notation:
$\P_\be=(\BZ[\be])^*$). The following claim follows from
\cite{FrolTa, Chapter~III, \bf{(2.20)}}.

\proclaim{Proposition 1.2} There exists $\xi_0\in\BQ(\be)\setminus\BZ[\be]$
such that
$$
\P_\be=\xi_0\cdot\BZ[\be].
$$
It can be given by the formula $\xi_0=(g'(\be))^{-1}$.
\endproclaim

\definition{Definition} Let the {\it Pisot group} $\A_\be$ be
defined as the quotient group $\P_\be/\BZ[\be]$.
\enddefinition

Let $\be_1=\be$, and $\be_2,\be_3,\dots,\be_m$
denote its conjugates. Finally, let
$D=D(\be)=\prod_{i\neq j}(\be_i-\be_j)^2$, i.e. the discriminant of $\be$.

\proclaim{Theorem 1.3} The order of the Pisot group is $|D(\be)|$.
\endproclaim
\demo{Proof} Again, we will use the facts from classical Number Theory.
By Dedekind's Ramification Theorem, for any lattice $M$ over $\BZ$, the
order of $M^*:M$ is $|D|$, where $D=\det M_\be$,
$$
M_\be=\pmatrix \Tr(1)         & \Tr(\be)   & \dots & \Tr(\be^{m-1}) \\
               \Tr(\be)       & \Tr(\be^2) & \dots & \Tr(\be^m)     \\
               \dots          & \dots      & \dots & \dots          \\
               \Tr(\be^{m-1}) & \Tr(\be^m) & \dots & \Tr(\be^{2m-1})
      \endpmatrix         \tag1.2
$$
(see \cite{FrolTa, Chapter~III, \bf{(2.8)}}).
It suffices to show that $\det M_\be=D$. Sometimes, this relation is
taken as a definition for the discriminant. Otherwise, observe that
$M_\be=V_\be\cdot V_\be^T$, where
$$
V_\be=\pmatrix 1           & 1           & \dots & 1          \\
               \be_1       & \be_2       & \dots & \be_m      \\
               \be_1^2     & \be_2^2     & \dots & \be_m^2    \\
               \dots       & \dots       & \dots & \dots      \\
               \be_1^{m-1} & \be_2^{m-1} & \dots & \be_m^{m-1}
      \endpmatrix,
$$
whence $\det M_\be=(\det V_\be)^2=\prod_{i\neq j}(\be_i-\be_j)^2=D$.\qed
\enddemo
The following simple fact answers the question about the finer structure
of the Pisor group. Let $d=\min\,\{l\ge1 :lM_\be^{-1}\in\BZ^{m\times m}\}$.

\proclaim{Lemma 1.4} The Pisot group $\A_\be$ contains $\BZ/d\BZ$ as a
subgroup. In particular, if $d=|D|$, then $\A_\be$ is cyclic, and
if $|D|/d$ is prime, then $\A_\be$ is isomorphic to
$(\BZ/d\BZ)\times(\BZ/(|D|/d)\BZ)$.
\endproclaim

\proclaim{Corollary 1.5} Assume the entries of $DM_\be^{-1}$ to be
coprime. Then $\A_\be$ is cyclic.
\endproclaim

\remark{Remark} All the above results are valid for {\it any} algebraic
unit, if one takes (1.1) as the definition of $\P_\be$.
\endremark

\example{Examples} {\bf 1.} Quadratic units. Let $\be^2=k\be\pm1$ with
$k\ge1$ for $+1$ and $k\ge3$ for $-1$. Then by direct inspection,
$$
M_\be=\pmatrix 2&k\\k&k^2\pm2 \endpmatrix,
$$
and $\xi_0=\frac1{\sqrt D}$. Then by Lemma~1.4, $\A_\be$ is isomorphic to
$\BZ/D\BZ$ if $k$ is odd, and to $(\BZ/(D/2)\BZ)\times(\BZ/2\BZ)$ otherwise.
The Pisot group for this case was considered in \cite{SV2}.
\smallskip
\noindent
{\bf 2.} Let $\be$ be the ``tribonacci number", i.e. the positive root
of $x^3=x^2+x+1$. Here $D=-44$,
$$
M_\be=\pmatrix 3&1&3\\1&3&7\\3&7&11\endpmatrix,\quad
M_\be^{-1}=\frac1{22}\pmatrix8&-5&1\\-5&-12&9\\1&9&-4\endpmatrix,
$$
$\xi_0^{-1}=g'(\be)=-1-2\be+3\be^2$, and
$\xi_0=\frac1{22}(1+9\be-4\be^2)$. By Lemma~1.4,
$\A_\be\cong(\BZ/22\BZ)\times(\BZ/2\BZ)$.
\smallskip
\noindent
{\bf 3.} Finally, let $\be$ be the smallest Pisot number, i.e. the root
of $\be^3=\be+1$ (see, e.g., \cite{DuPi}). Here $D=-23$,
$$
M_\be=\pmatrix 3&0&2\\0&2&3\\2&3&2 \endpmatrix,\quad
M_\be^{-1}=\frac1{23}\pmatrix5&-6&4\\-6&-2&9\\4&9&-6\endpmatrix,
$$
and since 23 is a prime, $\A_\be\cong \BZ/23\BZ$.
\endexample

At the end of the section, we would like to establish a link
bewteen the groups $\P_\be$ and $\A_\be$ and certain groups of recurrent
sequences, which will be used in the next section. Let
$$
R_\be=\{\{T_n\}_1^\infty\in\BZ^\infty\mid \exists j: T_{n+m}=k_1T_{n+m-1}+
k_2T_{n+m-2}+\dots+k_mT_n, \ n\ge j\}.
$$
Obviously, $R_\be$ is a group under addition.

\proclaim{Proposition 1.8} The groups $\P_\be$ and $R_\be$ are isomorphic.
\endproclaim
\demo{Proof} Let $\{T_n\}_1^\infty\in R_\be$. Put
$$
\xi:=\lim_{n\to+\infty}\be^{-n}T_n.
$$
Then $\xi\in \P_\be$, because $\|\xi\be^n\|=|\xi\be^n-T_n|\to0$ (in view
of $\be$ being a Pisot number, whence $|\xi-\be^nT_n|=o(\be^{-n}))$.

The inverse mapping from $\P_\be$ to $R_\be$ assigns to each $\xi\in\P_\be$
the sequence $\{T_n\}$, where $T_n$ is defined as the closest integer
for $\xi\be^n$. By the fact that $\be$ is a Pisot number, $T_n$ will
eventually satisfy the relation in question.\qed
\enddemo
Let now the equivalence relation $\sim$ on $R_\be$ be defined
as follows: $\{T_n\}\sim\{T_n'\}$ iff
$\lim_n \be^{-n}(T_n-T_n')\in\BZ[\be]$. Then obviously, the quotient
group $R_\be/\sim$ is isomorphic to $\A_\be$.

\head 2. Symbolic representation of the Pisot group
\endhead

\definition{Definition} A representation of an $x\in[0,1)$ of the form
$$
x=\sum_1^\infty\e_k\be^{-k}
$$
is called the {\it $\be$-expansion} of $x$ if the ``digits"
$\{\e_k\}_1^\infty$
are obtained by means of the greedy algorithm (similarly to the decimal
expanions), i.e., $\e_1=\e_1(x)=[\be x],\ \e_2=\e_2(x)=[\be\{\be x\}]$, etc.
The set of all possible sequences $\{\{\e_k(x)\}_1^\infty:x\in[0,1)\}$ is
called the (one-sided) {\it $\be$-compactum} and denoted by $X_\be^+$.
\enddefinition

The $\be$-compactum can be described more explicitly. Let
$1=\sum_1^\infty d_k\be^{-k}$ be the expansion of 1 defined by the greedy
algorithm, i.e.,
$d_1'=[\be],\ d'_2=\{\be[\be]\}$, etc. If the sequence $\{d_n'\}$
is not finite, we put $d_n\equiv d_n'$. Otherwise let $k=\max\,\{j:d_j'>0\}$,
and $(d_1, d_2, \dots):=(\ov{d_1',\dots, d'_{k-1},d'_k-1})$, where the bar
denotes a period.

We will write $\{x_n\}_1^\infty\prec\{y_n\}_1^\infty$
if $x_n<y_n$ for the least $n\ge1$ such that $x_n\neq y_n$. Then
by definition,
$$
X_\be^+=\{\{\e_n\}_1^\infty:(\e_n,\e_{n+1},\dots)\prec(d_1,d_2,\dots)
\ \text{for all}\ n\in\BN\}
$$
(see \cite{Pa}). Similarly, we define the {\it two-sided $\be$-compactum}
as
$$
X_\be=\{\{\e_n\}_{-\infty}^\infty:(\e_n,\e_{n+1},\dots)\prec(d_1,d_2,\dots)
\ \text{for all}\ n\in\BZ\}.
$$
Both compacta are naturally endowed with the weak topology, i.e. with the
topology of coordinate-wise convergence. For $\be$ Pisot the properties
of the $\be$-compactum are well-studied. Its main property is that
it is {\it sofic} (see, e.g., the review \cite{Bl}). Let $Fin(\be)$ denote
the set of $x$'s whose $\be$-expansions are {\it finite}, i.e. have the
tail $0^\infty$. Obviously, $Fin(\be)\subset\BZ[\be]\cap[0,1)$, but
the inverse inclusion does not hold for an arbitrary Pisot unit.

\definition{Definition} A Pisot unit is called {\it finitary}, if
$$
Fin(\be)=\BZ[\be]\cap[0,1).
$$
\enddefinition

\example{Examples} {\bf 1.} The Pisot number being the
principal root of the equation
$x^m=k_1x^{m-1}+\dots+k_m,\ k_1\ge k_2\ge\dots\ge k_m\ge1$,
is known to be finitary \cite{FrSo}.
\smallskip
\noindent
{\bf 2.} A quadratic Pisot unit is finitary if and only if
its norm is $+1$. Indeed, if $\be^2=k\be+1,\ k\ge1$, then the claim follows
from the previous one. If $\be^2=k\be-1,\ k\ge3$, then
in view of $(\be-1)^2=\be^2-2\be+1=(k-2)\be$, we have the $\be$-expansion
$\be-1=(k-2)\be^{-1}+(k-2)\be^{-2}+\dots$, whence
$\BZ[\be]\cap[0,1)\not\subset Fin(\be)$.
\smallskip
\noindent
{\bf 3.} For the cubic Pisot units also exists a full criterion
due to Sh.~Akiyama \cite{Ak}.
Namely, the norm of a finitary cubic $\be$ must be $+1$, i.e.
$\be^3=k_1\be^2+k_2\be+1$.
Then $\be$ is finitary if and only if $k_1\ge0$ and $-1\le k_2\le k_1+1$.
Thus, the tribonacci number and the smallest Pisot number
are both finitary, whereas, for example, the principal root of
$x^3=3x^2-2x+1$ is not.
\smallskip
\noindent
{\bf 4.} Finally, the second in order Pisot number, i.e. the positive root
of $x^4=x^3+1$, is non-finitary, as $\be^{-2}+\be^{-3}=\be^{-1}+\be^{-6}+
\be^{-11}+\dots$
\endexample

From here on we will assume $\be$ to be finitary.

\remark{Remark} The $\be$-expanions can be easily extended from
$[0,1)$ to $\BR_+$ in the way analogous to the decimal expansions.
This allows to add in $X_\be$ two sequences finite to the left.
Namely, if $\ov\e$ and $\ov\e'$ are both finite to the left,
put $x:=\sum_{-\infty}^\infty(\e_k+\e'_k)\be^{-k}$; then by definition,
$\ov\e+\ov\e'$ is the $\be$-expansion of $x$. The same is true
for subtraction, though $\ov\e-\ov\e'$ is well defined only
for the sequences, for which $\sum_k (\e_k-\e'_k)\be^{-k}\ge0$.
\endremark

\proclaim{Theorem} \rm{(A. Bertrand,  K. Schmidt)} For a Pisot $\be$,
any $x\in\BQ(\be)\cap\BR_+$ has an ultimately periodic $\be$-expansion.
\endproclaim

Therefore, all elements of $\P_\be\cap\BR_+$
have ultimately periodic $\be$-expanions.
Since their denominators in the standard basis of the field are bounded
(by $|D|$), these periods cannot be too long, whence it follows that
there is a finite set of such periods.

Our goal is to represent the Pisot group $\A_\be$ in $X_\be$. We perform
the following operation: to a $[\xi]\in\A_\be$ the mapping
$g$ assigns {\it all} purely periodic
two-sided sequences in $X_\be$ which are extensions to the left of the
periodic tails that occur in the $\be$-expansions of the nonnegative
$\xi'\sim\xi$. By the above argument, $\#g([\xi])<\infty$ for any
$\xi\in\P_\be\cap\BR_+$. Now we identify all the
images $\{g(\xi'): \xi'\in[\xi]\}$. This leads to
the mapping $g$ acting from $\A_\be$ onto the set that we will denote
by $\goth A_\be$.

\proclaim{Lemma 2.1} The mapping $g:\A_\be\to\goth A_\be$ is 1-to-1.
\endproclaim
\demo{Proof} We need to prove the injectivity of $g$ only. Let
$g([\xi])=g([\xi'])$. This means that there exist $\xi_1\in[\xi],\
\xi_2\in[\xi']$ whose tails of the $\be$-expansions coincide. But this
implies $\xi_1-\xi_2\in\BZ[\be]$, whence $[\xi]=[\xi']$.\qed
\enddemo

Thus, we have defined a certain finite arithmetic group which may be
regarded as a representation of the Pisot group in $X_\be$.
Although its elements are not sequences but certain
classes of equivalence, {\it arithmetically}
they represent one and the same sequence, because besides
arithmetic operations inherited from $\A_\be$, one may define
some internal arithmetic in $X_\be$. More precisely, addition
in $\goth A_\be$ can be carried out directly through the whole
compactum in the following way: let $\ov\e,\ov\e'\in X_\be$,
then $\ov\e+\ov\e'=\lim_{N\to+\infty}(\e^{(N)}+\e^{'(N)})$, where
$\e^{(N)}=(\dots,0,0,\dots,0,\e_{-N},\e_{-N+1},\dots)$, addition
of the sequences finite to the left being executed as explained in
the remark above.

Let us prove accurately that $\goth A_\be$ is a group under addition
in this sense as well. Let $Fin_k(\be)$ denote the set of sequences
$\ov\e$, for which $\e_j=0,\ j\ge k$. The theorem due to Ch.~Frougny
and B.~Solomyak says that if $\be$ is finitary, then there exists
$L=L(\be)$ such that $\ov\e, \ov\e'\in Fin_k(\be)$ implies
$\ov\e+\ov\e'\in Fin_{k+L}(\be)$, i.e. for the finitary Pisot numbers
the ``waiting time" is bounded \cite{FrSo}.

\proclaim{Lemma 2.2} Let $[\ov\e],[\ov\e']\in\goth A_\be$. Then regardless
of the choice of representatives in the equivalence classes,
$[\ov\e\pm\ov\e']$ is well defined.
\endproclaim
\demo{Proof} Let $\delta^{(N)}:=\e^{(N)}+\e^{'(N)}$. Then
$\delta^{(N+k)}-\delta^{(N)}\in Fin_{-N+L}(\be)$ by the theorem cited
above, whence the limit does exist. To obtain $-\ov\e$, we subtract
$\e^{(N)}$ from the sequence $(\dots,0,0,\dots,1,0,0,\dots)$, where 1
stands at the $(-N-1)$'th place and then pass to the limit. Even if it
did not exist, all partial limits would belong to one and the same
equivalence class, since the difference of their periodic parts
would belong to $\BZ[\be]$.\qed
\enddemo

Thus, $\goth A_\be$ is an additive group in $X_\be$ in the sense of its
natural arithmetics. Its obvious property is that it is
shift-invariant, i.e. contains any sequence together will all
its shifts.
Pursuing the idea from Section~1, we will give another way of obtaining
the sequences from $\goth A_\be$, which in a sense looks more traditional
(see \cite{FrSa1} for the case of the golden mean).

\proclaim{Lemma 2.3} For any element $[\ov\e]\in\goth A_\be$ there exists
an eventually nonnegative sequence $\{T_n\}_1^\infty\in
R_\be$ such that the set of partial limits for the sequence of
its $\be$-expansions $\{\ov\e_n\}_1^\infty$ is just the class
$[\e]$.
\endproclaim
\demo{Proof} Take $[\xi]=g^{-1}([\ov\e])$, and, as above, define $T_n$
as the closest integer to $\xi'\be^n$ for some nonnegative $\xi'\sim\xi$.
Then $T_n=\xi'\be^n+o(1)$, and $T_n$ are eventually nonnegative, as
$\xi'>0$. Hence it follows that to take the
set of partial limits of its $\be$-expansions is the same operation
as ``pulling out" the periodic tail for a $\xi\in\P_\be\cap\BR_+$.\qed
\enddemo

\example{Examples} {\bf 1.} For the quadratic Pisot units
$\be^2=k\be\pm1$ we proved in
Section~1 that $\P_\be=\frac1{\sqrt D}\cdot \BZ[\be]$. Since
$\frac1{\sqrt D}=\frac{k\be^2}{\be^4-1}$, all sequences in $\goth A_\be$
are of period~4. Let us consider some subcases. We will use the convention
to write just the period instead of the whole periodic sequence.
\smallskip
\noindent
{\bf 1.1.} $\be^2=\be+1$. Here $\goth A_\be=\{0000,1000,0100,0010,0001\}$.
Recall that by definition, $X_\be$ does not contain sequences ending by
$101010\dots$, that is why there is no need to identify
0000 with 1010 and 0101.
Thus, every class in $\goth A_\be$ consists of just one element.
If $F_1=1, F_2=2,\dots$ is the Fibonacci sequence, then by Lemma~2.3,
the set of partial limits for the $\be$-expansions of $(F_n)_1^\infty$
in $X_\be$ will give us just the four nonzero sequences
in $\goth A_\be$.\footnote{It is easy to see that the Fibonacci
sequence is the basis of the module $R_\be$, i.e. for any sequence
$\{T_n\}\in R_\be,\ T_n=\sum_{k=0}^s\e_kF_{n+k},\ \e_k\in\BZ,\ n\ge n_0$.}
However, this may be checked directly as well:
$$
\align
F_k&=\be^{k-1}+\be^{k-5}+\be^{k-9}+
\dots+\be^{-k+3}+\be^{-k},\quad\text{$k$ even},\\
F_k&=\be^{k-1}+\be^{k-5}+\be^{k-9}+
\dots+\be^{-k+5}+\be^{-k+1},\quad\text{$k$ odd},
\endalign
$$
whence the limits of $(F_{4k+j})_{k=1}^\infty$ in $X_\be$
do exist for $j=0,1,2,3$
and yield the four sequences in $\goth A_\be$. For more details see
\cite{FrSa1}, \cite{SV1}.
\smallskip
\noindent
{\bf 1.2.} $\be^2=2\be+1$. Here $\#\goth A_\be=8$, and
$\goth A_\be=\{0000,1010,0101,2000,0200,\allowmathbreak 0020,0002,1111\}$.
The elements 1010, 0101, 1111 are of order 2, whereas all other nonzero
elements are of order 4. A similar pattern takes place for
$\be^2=k\be+1,\ k \text{ even}$.

Note that in \cite{FrSa2} the authors study two similar groups for
arbitrary quadratic Pisot units. One of them, $H_\be$,
is the group of {\it all}
sequences of period 4; its order is $k^2D$, and it is in fact isomorphic
to $\frac{\BZ[\be]}{k\sqrt D}/\BZ[\be]$. Another one, $G_\be$, is, on the
contrary, smaller than $\goth A_\be$; it is related to a certain finite
automaton. For example, in the case $\be^2=2\be+1$, $G_\be=\{0000,0101,
1010,1111\}\cong (\BZ/2\BZ)\times(\BZ/2\BZ)$. In the general quadratic
case one can show that $G_\be\subset \goth A_\be\subset H_\be$, both
inclusions being, generally speaking, proper. However, in the case
of the golden mean (and only in this case), $G_\be=\goth A_\be=H_\be$.
\smallskip
\noindent
{\bf 2.} For the tribonacci number the situation with the periods
of the sequences in $\goth A_\be$ proves to be more complicated. Expanding
the elements of $\P_\be$, we see that these periods are 1,2,3 and 10.
More precisely, besides 0, $\goth A_\be$ consists of 40 sequences
of period 10, namely, of $\{1000110000,1010000110,1001011000,1001101100\}$
together with all their shifts, two sequences of period two: $\{01,10\}$
and {\it one} sequence of period 3: 100. One may ask: where are
all its shifts, mustn't they belong to $\goth A_\be$ as well? The
answer is simple: {\it arithmetcially} they are all equivalent, because
$(\be+\be^{-2}+\be^{-5}+\dots)-(1+\be^{-3}+\be^{-6}+\dots)=
\frac{\be-1}{1-\be^{-3}}=1\in\BZ[\be]$. This example shows that even for
a finitary $\be$
not necessarily any $\xi\in\P_\be$ and $\xi+l,\ l\in\BZ[\be]$ will have
one and the same tail. Therefore the definition of $\goth A_\be$ as
a quotient set is essential.
\smallskip
\noindent
{\bf 4.} $\be^3=\be+1$. Here $D=-23$.
Let us prove a simple claim about the structure of the group
$\goth A_\be$ in the case of prime $|D|$.

\proclaim{Lemma 2.4} If $|D|$ is a prime number,
then $|D|$ consists of $|D|-1$ sequences of period $|D|-1$ being the shifts
of each other, and 0.
\endproclaim
\demo{Proof} Since $\goth A_\be$ is shift-invariant and cyclic,
any $\ov\e\in\goth A_\be\setminus\{0\}$ belongs to it together with
all its shifts, and they are arithmetically non-equivalent.\qed
\enddemo
Return to the example. Here the period of $\xi_0$ is
$10000100000000100000000$, and $\goth A_\be$ consists of this sequence
together with its 21 shifts and 0.
\endexample

\head 3. Symbolic dynamics associated with the Pisot group
\endhead

Our goal is to show that a certain natural ``arithmetic" coding of
the ``companion" toral automorphism, is not
one-to-one a.e. and that its kernel is in fact coincides with the group
$\goth A_\be$. In this section we still assume $\be$ to be a
finitary Pisot unit.

We will need some facts from hyperbolic
dynamics. Let $T$ be a hyperbolic automorphism of the torus
$\BT=\BR^m/\BZ^m$, $L_s$ and $L_u$ denote respectively the leaves of the
stable and unstable foliations passing through $\0$.
Recall that a point homoclinic to $\0$ or simply a {\it homoclinic point}
is a point belonging to $L_s\cap L_u$. In other words, $\t$ is homoclinic iff
$T^n\t\to\0$ as $n\to\pm\infty$. The homoclinic points are a group
under addition isomorphic to $\BZ^m$, and we will denote it by
$\H(T)$. Each homoclinic point $\t$ can be obtained as follows:
take some $\bold n\in\BZ^m$ and project it onto $L_u$ along $L_v$
and then onto $\BT$ by taking the fractional parts of
all coordinates of the vector (see \cite{Ver}).

Let $T=T(\be)$ be the group automorphism of $\BT$ given
by the {\it companion matrix}, i.e. by
$$
M=M(\be)=\pmatrix k_1 & k_2 & k_3 & \dots & k_{m-1} & k_m \\
                  1   &  0  &  0  & \dots & 0       &  0  \\
                  0   &  1  &  0  & \dots & 0       &  0  \\
                  \dots & \dots & \dots& \dots & \dots &\dots \\
                  0   &  0  &  0 & \dots  & 1       & 0
         \endpmatrix.
$$
Then $(1,\be^{-1},\dots,\be^{-m+1})$ is an eigenvector corresponding
to the eigenvalue $\be$. In our case $L_u$ is one-dimensional,
whence $\t\in\H(T)$ iff $\t=(\xi,\xi\be^{-1},\dots,
\xi\be^{-m+1})\mod\BZ^m$ for some $\xi\in\P_\be$.

Consider two special homoclinic points. Let $\xi_0$ denote the
generator of the Pisot group defined in Proposition~1.4, and $\t_0$ denote
the fundamental homoclinic point given by the formula
$$
\t_0=(\xi_0,\xi_0\be^{-1},\dots,\xi_0\be^{-m+1})
$$
(we omit the natural projection of $\BR^m$ onto the torus, i.e.
write in the coordinates of $\BR^m$), and
$$
\t:=(1,\be^{-1},\dots,\be^{-m+1}).
$$
Consider the two mappings acting from $X_\be$ onto $\BT$, namely:
$$
\align
\phi_0(\e)&=\sum_{k\in\BZ}\e_kT^{-k}\t_0=
\lim_{N\to+\infty}\(\sum_{k=-N}^\infty \e_k\be^{-k}\)\pmatrix
\xi_0\\ \xi_0\be^{-1}\\ \dots\\ \xi_0\be^{-m+1}\endpmatrix
\mod\BZ^m, \\
\phi(\e)&=\sum_{k\in\BZ}\e_kT^{-k}\t=\lim_{N\to+\infty}\(\sum_{k=-N}^\infty
\e_k\be^{-k}\)\pmatrix 1\\ \be^{-1}\\ \dots\\ \be^{-m+1}\endpmatrix
\mod\BZ^m.
\endalign
$$

Both mappings are well defined, as $\|\be^N\|\to0,\ \|\xi_0\be^N\|\to0$
exponentially. They are obviously continuous, and their important property
is that they are bounded-to-one (see \cite{Sch} for $\phi_0$; in the case of
$\phi$ is essentially the same). Their main
value is that they both semiconjugate the shift $\tau$
on the compactum $X_\be$ and the automorphism $T$, i.e.,
$\phi_0\tau=T\phi_0,\ \phi\tau=T\phi$.
Thus, both mapping may be regarded as {\it arithmetic codings} of $T$.
The mapping $\phi_0$ was introduced in \cite{SV2} for $m=2$,
and in \cite{Sch}
for higher dimensions. At the same time, the mapping $\phi$ had been
historically first attempt of arithmetic coding of an automorphism of
the torus \cite{Ber}.

\proclaim{Theorem} {\rm (K. Schmidt \cite{Sch})} The mapping $\phi_0$ is
one-to-one on the set of doubly transitive sequences of $X_\be$, i.e.
on the sequences $\ov\e$ such that the sets $\{\tau^n\ov\e,\ n\ge k\}$
and $\{\tau^n\ov\e,\ n \le -k\}$ are both dense in $X_\be$ for every
$k\ge0$. Therefore, $\phi_0$ is bijective almost everywhere with respect to
the measure of maximal entropy for $\tau$.\footnote{The nature of this measure
is unimportant; in fact, one may take any shift-invariant measure that is
strictly positive on all cylinders $[\e_1=i_1,\dots,\e_k=i_k]\subset X_\be$.}
\endproclaim

Our goal in this section is to show that $\phi$ is $|D|$-to-1 a.e. and that
$\phi^{-1}(\0)$, after the identification of arithmetically equivalent
sequences, will coincide with $\goth A_\be$.

\definition{Definition}
We will say that $\ov\e$ is {\it equivalent} to $\ov\e'$ if
$\phi_0(\ov\e)=\phi_0(\ov\e')$.
\enddefinition

We will denote this equivalence relation by $\sim$. Let $X_\be'=X_\be/\sim$;
by Schmidt's theorem, $\#[\ov\e]=1$ for a.e. sequence $\ov\e$.
Let addition (subtraction) in $X_\be'$ be defined through
the torus, i.e., $[\ov\e]\pm[\ov\e']:=\phi_0^{-1}(\phi_0(\ov\e)
\pm\phi_0(\ov\e'))$. Obviously, $X_\be'$ is a group under addition isomorphic
to $\BT$.

\remark{Remark} In fact, to add $[\ov\e]$ to $[\ov\e']$, one may avoid
using the torus. Indeed, following the definition from the previous
section, consider the sequence
$\{\e^{(N)}+\e^{'(N)}\}_1^\infty$. It has got some set of partial limits
(possibly more than one), but they all belong to one and the same
equivalence class by the continuity of $\phi_0$. Thus, briefly,
$[[\ov\e]+[\ov\e']]=[\lim_N(\{\e^{(N)}\}+\{\e^{'(N)}\})]$.
\endremark

\proclaim{Theorem 3.1}
\roster
\item The mapping $\phi$ is well defined on the quotient compactum $X_\be'$.
\item $\phi: X_\be'\to\BT$ is a group homomorphism, its kernel being
equal to $\goth A_\be$.
\endroster
\endproclaim
\demo{Proof} (1) We need to show that if $\phi_0(\ov\e)=\phi_0(\ov\e')$, then
$\phi(\ov\e)=\phi(\ov\e')$. By the definition of $\phi_0$, we have
$\|\xi_0u_N\|\to0$, where $u_N=\sum_{k=-N}^\infty(\e_k-\e_k')\be^{-k}$.
By Proposition~1.2, $\xi_0^{-1}\in\BZ[\be]$, whence for {\it any}
$\{u_N\}$ such that $\|\xi_0u_N\|\to0$, $\|u_N\|\to0$ as well.
Hence by the definition of $\phi$, $\|\phi(\e^{(N)})-\phi(\e^{'(N)})\|\le
\text{const}\cdot\|u_N\|\to0$.\newline
(2) It is shown in \cite{Sch} that a mapping $\phi$ would be bounded-to-one
for any choice of homoclinic point $\t$. Since the set
$\Cal O=\phi^{-1}(\0)$ is finite and shift-invariant, it consists
of purely periodic sequences only. Let $\ov\e\in\Cal O$, its period
equal to $d$ and $\al=\sum_{k=1}^\infty\e_k\be^{-k}$. Then by the definition
of $\phi$, $\|\al\be^{dN}\|\to0$ as $N\to +\infty$. Considering the
sequences $\tau\ov\e,\ \tau^2\ov\e,\dots,\tau^{d-1}\ov\e$, we see
that $\|\al\be^{dN+j}\|\to0$ for $j=0,1,\dots,d-1$, whence
$\|\al\be^N\|\to0$ as $N\to+\infty$, i.e., $\al\in\P_\be$, and
the claim follows from Lemma~2.1. Conversely,
if $\al\in\P_\be$, then $\|\al\be^N\|\to0$, whence $\phi(\ov\e)=\0$.\qed
\enddemo

\proclaim{Corollary 3.2} The mapping $\phi$ is $|D|$-to-1 a.e.
\endproclaim

Let $A=\phi\phi_0^{-1}$. This mapping is well defined a.e. on the torus,
and we extend it to the whole $\BT$ by continuity. The following claim
is straightforward.

\proclaim{Lemma 3.3} The mapping $A:\BT\to\BT$ is an endomorphism of the
torus. Its determinant equals $\pm D$.
\endproclaim

In practice, to find $A$, one might expand $\xi_0^{-1}$ into the sum
$\sum_{k=0}^{m-1}a_k\be^k,\ a_k\in\BZ$. Then $A=\sum_{k=0}^{m-1}a_kT^k$.

\example{Examples} {\bf 1.} $\be^2=k\be\pm1$. Here, as we know,
$\xi_0=\frac1{\sqrt D},\ \xi_0^{-1}=\sqrt D=2\be-k$. Hence
$$
A=2T-kI=\pmatrix 2k&\pm2\\2&0\endpmatrix-\pmatrix k&0\\0&k\endpmatrix=
\pmatrix k&\pm2\\2&-k\endpmatrix,
$$
and $\det A=-D(\be)$.
\smallskip
\noindent
{\bf 2.} $\be^3=\be^2+\be+1$. Here $\xi_0^{-1}=-1-2\be+3\be^2$, whence
$$
A=-I-2\pmatrix 1&1&1\\1&0&0\\0&1&0\endpmatrix+
3\pmatrix 1&1&1\\1&0&0\\0&1&0\endpmatrix^2=\pmatrix 3&4&1\\1&2&3\\3&-2&-1
\endpmatrix,
$$
and $\det A=44=-D(\be)$.

\endexample

\subhead Acknowledgement \endsubhead I wish to thank
Anatoly Vershik for our fruitful collaboration in dimension two
that has given rise to many ideas in the present paper for higher
dimensions. An essential part of this
work was written in June-July 1999 during my stay in LIAFA at
University Paris-VII, and I am greatly indebted
to Christiane Frougny for the invitation and stimulating
discussions. I am grateful to Shigeki Akiyama for his professional
number-theoretic advices and his active interest in this paper.

\enddocument